\documentclass[11pt]{article}
\usepackage{latexsym}
\usepackage{amsmath}
\begin{document}
 \setlength{\baselineskip}{18pt plus2pt minus2pt}

 \begin{center}

{\bf {\Large On bicyclic graphs with maximal Graovac-Ghorbani index}}\footnote{This work is supported by the National Natural Science Foundation of China (grant number 12061074, 12201471), the Doctoral Research Foundation of Longdong University (XYBYZK2306), the Special Foundation in Key Fields for Universities of Guangdong Province (2022ZDZX1034).}

\vskip 0.5cm

Rui Song$^{a,b}$, {Saihua Liu$^c$}\footnote{Corresponding author.}, Jianping Ou$^c$ \\

\vskip 0.2cm

{\small $^a$School of Mathematics and Information Engineering, Longdong University, Qingyang, Gansu 745000, P.R.China\\
$^b$Institute of Applied Mathematics, Longdong University, Qingyang, Gansu 745000, P.R.China\\
$^c$Department of Mathematics, Wuyi University, Jiangmen, Guangdong 529020, P.R.China\\
}

E-mail: 15819931748@163.com, lsh1808@163.com, oujp@263.net

\end{center}

{\bf Abstract}\ Graovac-Ghorbani index is a new version of the atom-bond connectivity index.
D. Pacheco et al. [D. Pacheco, L. de Lima, C. S. Oliveira, On the Graovac-Ghorbani Index for Bicyclic Graph with No Pendent Vertices, MATCH Commun. Math. Comput. Chem. 86 (2021) 429-448] conjectured a sharp lower and
upper bounds to the Graovac-Ghorbani index for all bicyclic graphs. Motivated by their nice work, in this paper we determine the maximal Graovac-Ghorbani index of bicyclic graphs and characterize the corresponding extremal graphs, which solves one of their Conjectures.


\section{ Introduction}

Molecular descriptors play a significant role in mathematical chemistry especially in the QSPR/QSAR investigations.
Among them, special place is revised for so-called topological indices [2], where topological index are numbers
associated with chemical structures as a tool for compact and effective description of structural formulas
used to study and predict the structure-property correlation of organic compounds[3-6]. It is known that
connectivity index has a widely application [7]. In [8], Estrada et al proposed
the concept of atom-bond connectivity index (ABC index) of a simple undirected graph $G$  as
$$
ABC(G)=\sum\limits_{uv\in E(G)}\sqrt{\frac{d(u)+d(v)-2}{d(u)d(v)}},
$$
where $E(G)$ is the edge set of graph $G$ and $d(u)$ is the degree of vertex $u$. ABC index has proven
to be a valuable predictive index in study of the heat of formation in alkanes[8-9]. And many
mathematical properties of this index was reported in [10-18].

Recently, Graovac and Ghorbani, defined a new version of the atom-bond connectivity index in [1] as follows:
$$
ABC_{GG}(G)=\sum\limits_{uv\in E(G)}\sqrt{\frac{n_u+n_v-2}{n_un_v}},
$$
where the summation goes over all edges in  graph $G$, $n_u$ denotes the number of vertices of $G$
whose distances to vertex $u$ are smaller than those to other vertex $v$ of the edge $e=uv$,
and $n_v$ defines similarly. This index is also called the second atom-bond connectivity index
of graph $G$ and is labeled as $ABC_2(G)$ [1,19-22]. Boris Furtula  think that which need be
denoted as $ABC_{GG}$ and be called as Graovac-Ghorbani index because of their essential
difference [23]. Therefore, this index is denoted as $ABC_{GG}$ and is called Graovac-Ghorbani
index in this paper. In the chemistry applications, $ABC_{GG}$ is used to model both the
boiling point and melting point of the molecules. Hence, it is also applied to the pharmaceutical field.
Some papers have given a lot of properties and have calculated many upper bounds for the second
atom-bond connectivity index [1, 19-23]. 
In [24], D. Pacheco, L. D. Lima and C. S. Oliveira determined the extremal graph with minimum Graovac-Ghorbani index among all bicyclic graphs with no pendent vertices. Additionally, they conjectured a sharp lower and upper bounds to the Graovac-Ghorbani index for all bicyclic graphs.
In this paper, we determine the maximal $ABC_{GG}$ index
of bicyclic graphs and characterize the extremal graphs.

\section{Preliminaries}

This section will present some lemmas. For two real numbers $x,y \geq 1$ and integer $a \geq 2$,
let $f(x,y)=\sqrt{\frac{x+y-2}{xy}}$, $g_a(x)=f(x,a)-f(x,a-1)$. Then
\begin{eqnarray*}
g_a (x) &=& f(x,a)-f(x,a-1) \nonumber \\
&=& \frac{2-x}{\sqrt{a(a-1)x}(\sqrt{(a-1)(x+a-2)}+\sqrt{a(x+a-3)})}.
\end{eqnarray*}
In [16] the authors show the following lemma. 

{\bf Lemma 2.1}[16]. \ $f(x,1)=\sqrt{\frac{x-1}{x}}$ is strictly increasing for $x$, $f(x,2)=\sqrt{\frac{1}{2}}$,
and $f(x,y)$ is strictly decreasing for $x$ if $y\geq3$.

By Lemma 2.1, we have the following Lemma 2.2, part of which has been given in [16,17],
so we omite its proof herein.

{\bf Lemma 2.2}[16,17].\ If $a \ge 2$, then $g_a(x)\leq g_a(2) \leq g_a(1)$ for $x \geq2$;
$g_a(1)$ is strictly decreasing for $a$, and $g_a (x)$ is strictly increasing for $a$ if $x \ge 3$.

\section{Bicyclic graphs with maximum $ABC_{GG}$ index}

Let $B_{n,p}$ be the set of bicyclic graphs with $n$ vertices and $p$ pendent vertices, clearly we have that
$0 \leq p \leq n-4$. Let $S(m_1,\cdots,m_k)$ be the unicyclic graph with cycle $v_1v_2\cdots v_kv_1$
and $m_i$ pendent vertices adjacent to vertex $v_i$ for all $i=1,2,\cdots,k$. Bicyclic graph
$S_n^{r,t} (m_1,\cdots,m_{r-1},n_1,\cdots,n_{t-1},m_0 )$ is
defined in following Figure 1  and $S_n^{3,3} (m_1,n_1,m_0 )=S_n^{3,3}(m_1,0,n_1,0,m_0)$.

\setlength{\unitlength}{1cm}
\begin{center}
\begin{picture}(12,6)
\put(6,0){Figure \ 1}
\put(0, 0.8){$S_n^{r,t}(m_1,m_2,...,m_{r-1},n_1,n_2,...,n_{t-1},m_0)$}
\put(3,3.5){\oval(1,1)[r]}
\put(4.5,3.5){\oval(2,2)[l]}
\bezier{10}(3,4)(2,3.5)(3,3)
\bezier{20}(4.5,4.5)(6.5,3.5)(4.5,2.5)
\put(3,3){\circle*{0.1}}
\put(3,4){\circle*{0.1}}
\put(3.5,3.5){\circle*{0.1}}
\put(4.5,4.5){\circle*{0.1}}
\put(4.5,2.5){\circle*{0.1}}
\multiput(4,4)(0,-0.5){3}{\circle*{0.1}}
\put(3.5,3.5){\line(1,1){0.5}}
\put(3.5,3.5){\line(1,0){0.5}}
\put(3.5,3.5){\line(1,-1){0.5}}
\multiput(4.5,5)(0.5,0){2}{\circle*{0.1}}
\multiput(4.5,2)(0.5,0){2}{\circle*{0.1}}
\multiput(4.5,4.5)(0,-2.5){2}{\line(0,1){0.5}}
\put(4.5,4.5){\line(1,1){0.5}}
\put(4.5,2.5){\line(1,-1){0.5}}
\put(4.4,5){$\overbrace{\dots}\limits^{n_1}$}
\put(4.4,2){$\underbrace{\dots}\limits_{n_{t-1}}$}
\multiput(2.5,2.5)(0.5,0){2}{\circle*{0.1}}
\multiput(2.5,4.5)(0.5,0){2}{\circle*{0.1}}
\multiput(3,2.5)(0,1.5){2}{\line(0,1){0.5}}
\put(3,3){\line(-1,-1){0.5}}
\put(3,4){\line(-1,1){0.5}}
\put(2.4,2.5){$\underbrace{\dots}\limits_{m_{r-1}}$}
\put(2.4,4.5){$\overbrace{\dots}\limits^{m_1}$}
\small{
\put(4.1,4){$1$}
\put(4.1,3.4){$2$}
\put(4,3.1){\vdots}
\put(4.1,2.9){$m_0$}
}
\put(2.7,3.5){$C_r$}
\put(5,3.5){$C_t$}

\multiput(8,5)(2,0){2}{\circle{1}}
\multiput(8.5,5)(1,0){2}{\circle*{0.1}}
\put(8.5,5){\line(1,0){1}}
\put(8.6,5.1){$u$}
\put(9.3,5.1){$v$}
\multiput(8.5,2.5)(1,0){2}{\circle{1}}
\put(9,2.5){\circle*{0.1}}
\put(9,3.9){$\downarrow$}
\put(7,0.8){Edge-lifting transformation on $uv$}
\put(7.8,5){$G_1$}
\put(9.8,5){$G_2$}
\put(8.3,2.5){$G_1$}
\put(9.3,2.5){$G_2$}
\put(9,3){\circle*{0.1}}
\put(9,2.5){\line(0,1){0.5}}
\end{picture}
\end{center}

Part of the following Lemma 3.1 has been shown in [19].

{\bf Lemma 3.1.}  Let $x, n$ be two positive integers with $1\leq x \leq n-6$.

	(1) If $n=7$, $x=1$, then
	$$	\sqrt{\frac{x}{x+1}}+\sqrt{\frac{n-3-x}{n-2-x}}+\sqrt{\frac{n-3}{(x+1)(n-2-x)}}=\sqrt{\frac{3}{4}}+\sqrt{2}.
	$$

	(2) If $n=8$ , then
$$
\sqrt{\frac{x}{x+1}}+\sqrt{\frac{n-3-x}{n-2-x}}+\sqrt{\frac{n-3}{(x+1)(n-2-x)}} \leq \sqrt{\frac{2}{3}}+
\sqrt{\frac{3}{4}}+\sqrt{\frac{5}{12}},
$$
with equality holding if and only if $x=2$.

	(3) If $n=9$ , then
$$
\sqrt{\frac{x}{x+1}}+\sqrt{\frac{n-3-x}{n-2-x}}+\sqrt{\frac{n-3}{(x+1)(n-2-x)}} \leq \sqrt{3}+
\sqrt{\frac{3}{8}},
$$
with equality holding if and only if $x=3$.

	(4) If $10 \leq n \leq 15$ , then
$$
\sqrt{\frac{x}{x+1}}+\sqrt{\frac{n-3-x}{n-2-x}}+\sqrt{\frac{n-3}{(x+1)(n-2-x)}} \leq \sqrt{\frac{2}{3}}+
\sqrt{\frac{n-5}{n-4}}+\sqrt{\frac{n-3}{3(n-4)}},
$$
with equality holding if and only if $x=2$.

	(5) If $n \ge 16 $ , then
$$
\sqrt{\frac{x}{x+1}}+\sqrt{\frac{n-3-x}{n-2-x}}+\sqrt{\frac{n-3}{(x+1)(n-2-x)}} \leq \sqrt{2}+
\sqrt{\frac{n-4}{n-3}}
$$
with equality holding if and only if $x=1$.

{\bf Proof.}  Let $f(n, x)= \sqrt{\frac{x}{x+1}}+\sqrt{\frac{n-3-x}{n-2-x}}+\sqrt{\frac{n-3}{(x+1)(n-2-x)}}$.
Statement (1) follows by direct calculation. If $n=8$ then $x=1$ or 2, (2) follows by comparing $f(8,1)$
and $f(8,2)$.  If  $n=9$ then $x=1, 2$ or 3, (3) also follows by comparing $f(9,1)$, $f(9,2)$ and $f(9,3)$.
So, we need only prove other parts in what follows.

Notice that $f(n,x)=f(n,n-3-x)$ for $x \in [1,n-4]$. In [19, Theorem 2.4], the authors show that
for $\frac{n-3}{2}\leq x \leq n-4$, if $10 \leq n \leq 15$ then $f(n,x)\leq  \sqrt{\frac{2}{3}}+
\sqrt{\frac{n-5}{n-4}}+\sqrt{\frac{n-3}{3(n-4)}}=f(n,n-5)$; if $n \geq 16$ then $f(n,x)\leq \sqrt{2}+
\sqrt{\frac{n-4}{n-3}}=f(n,n-4)$. The above discussion shows that for all $x$ with $1\leq x \leq n-6$,
if $10\leq n \leq 15$ then $f(n,x)\leq  \sqrt{\frac{2}{3}}+
\sqrt{\frac{n-5}{n-4}}+\sqrt{\frac{n-3}{3(n-4)}}=f(n,n-5)=f(n,2)$; if $n \geq 16$ then $f(n,x)\leq \sqrt{2}+
\sqrt{\frac{n-4}{n-3}}=f(n,n-4)=f(n,1)$. Hence, (4) and (5) follows. $\Box$

{\bf Edge-lifting transformation on edge $uv$ of graph $G$} [19]. Let $uv$ be a cut edge of a connected graph
$G$ but $uv$ be not a pendent edge. We delete edge $uv$ from $G$ at first, then identify vertices $u$
and $v$, and finally attach a new isolated vertex to this identified vertex to obtain a new graph. This
graph transformation is called an edge-lifting transformation on edge $uv$ of graph $G$.  This transformation
is pictured in Figure 1. 

In [19], the authors prove the following lemma for the case when $G$ is a unicyclic graph.
Here, we generalize it to connected graphs. 

{\bf Lemma 3.2. } If $G^\prime$ is the graph obtained by performing edge-lifting transformation
on edge $xy$ of graph $G$, then $$ABC_{GG}(G) < ABC_{GG}(G^\prime).$$

{\bf proof.}\  If $uv \in E(G_1) \cup E(G_2)$ (refer to Figure 1), by the definition of  $ABC_{GG}$ index
we have that
$$
\sqrt{\frac{n_u+n_v-2}{n_un_v}}
$$
contributes the same to $ABC_{GG}(G)$ and $ABC_{GG}(G^\prime)$.
Since $|G_1|\ge 2, |G_2| \ge 2$ and $|G_1|+|G_2|=|G|$, we have
$$
ABC_{GG}(G)-ABC_{GG}(G^\prime)=\sqrt{\frac{|G_1|+|G_2|-2}{|G_1||G_2|}}-\sqrt{\frac{|G|-2}{(|G|-1)\cdot1}}<0.
$$
And so, $ABC_{GG}(G)<ABC_{GG}(G^\prime)$.  $\Box$

Let $S_n^{r,t}$ be the set of such bicyclic graphs of order $n$ that have two cycles with length $r$
and $t$ respectively, furthermore, these two cycles have unique common vertex. 

{\bf Lemma 3.3. } If $G \in S_n^{r,t}$ is a connected bicyclic graph of order $n \ge 7$ then
$$
ABC_{GG}(G)\leq (n-5)\sqrt{\frac{n-2}{n-1}}+6\sqrt{\frac{1}{2}}
$$
with the equality holding only if $G= S_n^{3,3}(m_1,n_1,n-m_1-n_1-5)$ for some integers
$m_1,m_2,n_1,n_2,m_0 \ge 0$.

{\bf Proof.}\ If $G$ is not isomorphic to any $S_n^{r,t} (m_1,\cdots,m_{r-1},n_1,\cdots,n_{t-1},m_0 )$,
then $G$ has a non-pendent edge $uv$ that $uv \notin E(C_r ) \cup E(C_t)$, where $C_r$ and $C_t$
are the two cycles of graph $G$. By performing
edge-lifting transformation on edge $uv$ of $G$ we obtained a new graph $G_1$, from Lemma 3.2 it
follows that $ABC_{GG}(G) < ABC_{GG}(G_1)$. So, we may assume in what follows that
$G=S_n^{r,t}(m_1,\cdots,m_{r-1},n_1,\cdots,n_{t-1},m_0 )$ for some integers $r,t \geq3$. Now,  three
different cases occur.

Case 1. $r>3$, $t>3$. In this case, $G$ has $n -r-t+1 \le n-7$ pendent edges
and at least $r+t \ge 8$ non-pendent edges. Since $r>3$, $t>3$, for each non-pendent edge
$uv \in E(C_r)\cup E(C_t)$ we have $n_u \geq 2$ and $n_v \geq 2$.  by Lemma 2.1,
$\sqrt{\frac{n_u+n_v-2}{n_u n_v }}\leq \sqrt{\frac{1}{2}}<\sqrt{\frac{n-2}{n-1}}$  since $n \ge 7$
in this case. For any pendent edge $uv \in E(G)$, $\sqrt{\frac{n_u+n_v-2}{n_u n_v }}=
\sqrt{\frac{n-2}{n-1}}$. Hence,
\begin{eqnarray}
ABC_{GG}(G)&=&\sum\limits_{\substack{uv\in E(G) \\ d_u=1}}\sqrt{\frac{n_u+n_v-2}{n_un_v }}+
\sum\limits_{\substack {uv \in E(G) \\ d_{u},d_{v}\neq1}}\sqrt{\frac{n_u+n_v-2}{n_u n_v }}\nonumber \\
& \le & (n-r-t+1)\sqrt{\frac{n-2}{n-1}}+(r+t)\sqrt{\frac{1}{2 }}\nonumber \\
& = &(n-5)\sqrt{\frac{n-2}{n-1}}+(6-r-t)\sqrt{\frac{n-2}{n-1}}+(r+t)\sqrt{\frac{1}{2}}\nonumber \\
&<& (n-5)\sqrt{\frac{n-2}{n-1}}+6\sqrt{\frac{1}{2}}. \nonumber
\end{eqnarray}

Case 2. $r>3$, $t=3$. If $n_1,n_2\geq 1$, then the above inequality is also true; If $n_1,n_2=0$,
since $r \geq4$ it follows that $(r-6)\sqrt{\frac{1}{2}}-(r-5)\sqrt{\frac{n-2}{n-1}}<0$. Recalling
$\sqrt{\frac{1}{2}}<\sqrt{\frac{n-2}{n-1}}$, we deduce that
\begin{eqnarray}
ABC_{GG}(G)&=&\sum\limits_{uv\in E(C_r)}\sqrt{\frac{n_u+n_v-2}{n_un_v }}
+2\sqrt{\frac{n-3}{n-2}}+(n-r-2)\sqrt{\frac{n-2}{n-1}}\nonumber \\
&\leq&r\sqrt{\frac{1}{2}}+2\sqrt{\frac{n-3}{n-2}}-(r-3)\sqrt{\frac{n-2}{n-1}}+(n-5)\sqrt{\frac{n-2}{n-1}}\nonumber \\
&<&6\sqrt{\frac{1}{2}}+(r-6)\sqrt{\frac{1}{2}}-(r-5)\sqrt{\frac{n-2}{n-1}}+(n-5)\sqrt{\frac{n-2}{n-1}}\nonumber \\
&<& (n-5)\sqrt{\frac{n-2}{n-1}}+6\sqrt{\frac{1}{2}}. \nonumber
\end{eqnarray}

Finally, if $n_1\geq1, n_2=0$, since $n_1 \le n-7$ in this case it follows that
$$
\sqrt{\frac{n-3}{(n_1+1)(n-n_1-2)}}\leq \sqrt{\frac{1}{2}}\leq \sqrt{\frac{n_1}{n_1+1}}<\sqrt{\frac{n-2}{n-1}},
$$
and
$$
\sqrt{\frac{n-n_1-3}{n-n_1-2}}<\sqrt{\frac{n-2}{n-1}}.
$$
When $r \geq 5$ we have
\begin{eqnarray*}
ABC_{GG}(G)&=&\sum\limits_{uv\in E(C_r)}\sqrt{\frac{n_u+n_v-2}{n_un_v }}+\sqrt{\frac{n-3}{(n_1+1)(n-n_1-2)}} \\
&&+\sqrt{\frac{n_1}{n_1+1}}+\sqrt{\frac{n-n_1-3}{n-n_1-2}}+(n-r-2)\sqrt{\frac{n-2}{n-1}} \\
&<&r\sqrt{\frac{1}{2}}+\sqrt{\frac{1}{2}}+2\sqrt{\frac{n-2}{n-1}}+(n-r-2)\sqrt{\frac{n-2}{n-1}} \\
&=&6\sqrt{\frac{1}{2}}+(r-5)\sqrt{\frac{1}{2}}-(r-5)\sqrt{\frac{n-2}{n-1}}+(n-5)\sqrt{\frac{n-2}{n-1}} \\
&\leq& (n-5)\sqrt{\frac{n-2}{n-1}}+6\sqrt{\frac{1}{2}}.
\end{eqnarray*}

When $r=4$, let $f(x)=\sqrt{\frac{x}{x+1}}+\sqrt{\frac{n-3-x}{n-2-x}}+\sqrt{\frac{n-3}{(x+1)(n-2-x)}}$,
where $1\leq x \leq n-6$ and  $n \geq7$.  By Lemma 3.1, if $n=7$, $x=1$ then $f(1)=\sqrt{2}+\sqrt{\frac{3}{4}}
<\sqrt{2}+\sqrt{\frac{5}{6}}$; if $n=8$ then $f(x)\leq f(2)=\sqrt{\frac{2}{3}}+\sqrt{\frac{3}{4}}+
\sqrt{\frac{5}{12}}\approx2.3280<\sqrt{2}+\sqrt{\frac{6}{7}}\approx2.3400$;
if $n=9$ then $f(x)\leq f(3)=\sqrt{3}+\sqrt{\frac{3}{8}}\approx2.3444<\sqrt{2}+\sqrt{\frac{7}{8}}\approx2.3496$;
if $10 \leq n \leq 15$ then $f(x)\leq f(2)=\sqrt{\frac{2}{3}}+\sqrt{\frac{n-5}{n-4}}+\sqrt{\frac{n-3}{3(n-4)}}$.
Let $g(n)=\sqrt{\frac{2}{3}}+\sqrt{\frac{n-5}{n-4}}+\sqrt{\frac{n-3}{3(n-4)}}-\sqrt{2}-\sqrt{\frac{n-2}{n-1}}$.
Then $g(10)=\sqrt{\frac{2}{3}}+\sqrt{\frac{5}{6}}+\sqrt{\frac{7}{18}}-\sqrt{2}-\sqrt{\frac{8}{9}}\approx
-0.0040<0$, $g(11)=\sqrt{\frac{2}{3}}+\sqrt{\frac{6}{7}}+\sqrt{\frac{8}{21}}-\sqrt{2}-\sqrt{\frac{9}{10}}
\approx-0.0034<0$, $g(12)=\sqrt{\frac{2}{3}}+\sqrt{\frac{7}{8}}+\sqrt{\frac{9}{24}}-\sqrt{2}-
\sqrt{\frac{10}{11}}\approx-0.0034<0$, $g(13)=\sqrt{\frac{2}{3}}+\sqrt{\frac{8}{9}}+\sqrt{\frac{10}{27}}
-\sqrt{2}-\sqrt{\frac{11}{12}}\approx-0.0038<0$, $g(14)=\sqrt{\frac{2}{3}}+\sqrt{\frac{9}{10}}+
\sqrt{\frac{11}{30}}-\sqrt{2}-\sqrt{\frac{12}{13}}\approx-0.0043<0$, $g(15)=\sqrt{\frac{2}{3}}+
\sqrt{\frac{10}{11}}+\sqrt{\frac{12}{33}}-\sqrt{2}-\sqrt{\frac{13}{14}}\approx-0.0049<0$;
if $n \geq 16$ then $f(x)\leq f(1)=\sqrt{2}+\sqrt{\frac{n-4}{n-3}}<\sqrt{2}+\sqrt{\frac{n-2}{n-1}}$.
Hence, $f(x)<\sqrt{2}+\sqrt{\frac{n-2}{n-1}}$ for all $x$ with $1\leq x \leq n-6$. So,
\begin{eqnarray*}
ABC_{GG}(G)&=&\sum\limits_{uv\in E(C_r)}\sqrt{\frac{n_u+n_v-2}{n_un_v }}+\sqrt{\frac{n-3}{(n_1+1)(n-n_1-2)}} \\
&&+\sqrt{\frac{n_1}{n_1+1}}+\sqrt{\frac{n-n_1-3}{n-n_1-2}}+(n-r-2)\sqrt{\frac{n-2}{n-1}} \\
&\leq&6\sqrt{\frac{1}{2}}+(n-5)\sqrt{\frac{n-2}{n-1}}+\sqrt{\frac{n-3}{(n_1+1)(n-n_1-2)}} \\
&&+\sqrt{\frac{n_1}{n_1+1}}+\sqrt{\frac{n-n_1-3}{n-n_1-2}}-2\sqrt{\frac{1}{2}}-\sqrt{\frac{n-2}{n-1}} \\
&<& (n-5)\sqrt{\frac{n-2}{n-1}}+6\sqrt{\frac{1}{2}}.
\end{eqnarray*}
The proof is finished.   $\Box$

{\bf Lemma 3.4.} If $G \in S_n^{3,3}(m_1,m_2,n_1,n_2,m_0)$ has maximal $ABC_{GG}$ index, then
$\min \{m_1,m_2\}=0$ and $\min \{n_1,n_2\}=0$. 

{\bf Proof.} If $n \le 6$, the lemma is obviously true. So, in what follows we assume that $n \ge 7$. Firstly,
we deduce that if $m_1, n_1 \ge 1$ then
\begin{eqnarray}
&& ABC_{GG}(S_n^{3,3}(m_1,n_1,m_0))= \nonumber \\
& &\left(\sqrt{\frac{m_1}{m_1+1}}+\sqrt{\frac{n-m_1-3}{n-m_1-2}}+
\sqrt{\frac{n-3}{(m_1+1)(n-m_1-2)}}\right)\nonumber \\
& &+\left(\sqrt{\frac{n_1}{n_1+1}}+\sqrt{\frac{n-n_1-3}{n-n_1-2}}+
\sqrt{\frac{n-3}{(n_1+1)(n-n_1-2)}}\right)\nonumber \\
& &+(n-5)\sqrt{\frac{n-2}{n-1}} \\
& &>(n-5)\sqrt{\frac{n-2}{n-1}}+6\sqrt{\frac{1}{2}}. \nonumber
\end{eqnarray}

When  $m_1,m_2,n_1,n_2>0$, for each non-pendent edge $uv \in E(G)$ we have $n_u \geq2$ and $n_v \geq2$,
$\sqrt{\frac{n_u+n_v-2}{n_u n_v }}\leq \sqrt{\frac{1}{2}}$. Hence,
\begin{eqnarray}
ABC_{GG}(G)&=&ABC_{GG}(S_n^{3,3}(m_1,m_2,n_1,n_2,m_0))\nonumber \\
&=&\sum\limits_{\substack {uv\in E(G) \\ d_{u}=1}}\sqrt{\frac{n_u+n_v-2}{n_un_v }}+
\sum\limits_{\substack {uv\in E(G) \\ d_{u},d_{v}\neq1}}\sqrt{\frac{n_u+n_v-2}{n_u n_v }}\nonumber \\
& \leq & (n-5)\sqrt{\frac{n-2}{n-1}}+6\sqrt{\frac{1}{2}}. \nonumber
\end{eqnarray}

When exactly one of $m_1,m_2,n_1,n_2$ equals zero, say, $m_2=0$, by Lemma 2.1 we have

 $\sqrt{\frac{n_1+n_2}{(n_1+1)(n_2+1)}}\leq \sqrt{\frac{1}{2}}\leq
\sqrt{\frac{n_1}{n_1+1}}$,

$ \sqrt{\frac{n-n_1-3}{(n_2+1)(n-n_1-n_2-2)}}\leq \sqrt{\frac{1}{2}}\leq
\sqrt{\frac{n-n_1-3}{n-n_1-2}}$,

$\sqrt{\frac{n-n_2-3}{(n_1+1)(n-n_1-n_2-2)}}<\sqrt{\frac{n-3}{(n_1+1)(n-n_1-2)}}. $ \\
Hence,
\begin{eqnarray}
ABC_{GG}(G)&=&ABC_{GG}(S_n^{3,3}(m_1,0,n_1,n_2,m_0))\nonumber \\
&=& \sqrt{\frac{m_1}{m_1+1}}+\sqrt{\frac{n-m_1-3}{n-m_1-2}}+\sqrt{\frac{n-3}{(m_1+1)(n-m_1-2)}} \nonumber\\
& &+\sqrt{\frac{n_1+n_2}{(n_1+1)(n_2+1)}}+\sqrt{\frac{n-n_2-3}{(n_1+1)(n-n_1-n_2-2)}} \nonumber\\
& &+\sqrt{\frac{n-n_1-3}{(n_2+1)(n-n_1-n_2-2)}}+(n-5)\sqrt{\frac{n-2}{n-1}}\nonumber\\
&< & \left(\sqrt{\frac{m_1}{m_1+1}}+\sqrt{\frac{n-m_1-3}{n-m_1-2}}+\sqrt{\frac{n-3}{(m_1+1)(n-m_1-2)}}\right) \nonumber\\
& &+\left(\sqrt{\frac{n_1}{n_1+1}}+\sqrt{\frac{n-n_1-3}{n-n_1-2}}+\sqrt{\frac{n-3}{(n_1+1)(n-n_1-2)}}\right)\nonumber\\
& &+(n-5)\sqrt{\frac{n-2}{n-1}}\nonumber \\
&=&ABC_{GG}(S_n^{3,3}(m_1,n_1,m_0)). \nonumber
\end{eqnarray}

When exactly two of $m_1,m_2,n_1,n_2$ equal zero, if $m_1=m_2=0$ (the case when $n_1=n_2=0$ is similar)
then $G$ has unique non-pendent edge $uv$ such that $n_u=n_v=1$, which contributes 0 to the index;
$G$ also has two such non-pendent edges $uv$ and $xy$ that $n_u=1$ and $n_v=n-2$,
each of them contributes $\sqrt{\frac{n_u+n_v-2}{n_un_v }}=\sqrt{\frac{n-3}{n-2}}$
to the index. Since $2\sqrt{\frac{n-3}{n-2}}-3\sqrt{\frac{1}{2}}=2(\sqrt{\frac{n-3}{n-2}}-\sqrt{\frac{9}{8}})<0$,
it follows that $2\sqrt{\frac{n-3}{n-2}}<3\sqrt{\frac{1}{2}}$. So,
\begin{eqnarray}
ABC_{GG}(G)&=&ABC_{GG}(S_n^{3,3}(0,0,n_1,n_2,m_0))\nonumber \\
&=&\sum\limits_{\substack {uv\in E(G) \\ d_{u}=1}}\sqrt{\frac{n_u+n_v-2}{n_un_v }}+
\sum\limits_{\substack{ uv\in E(G) \\ d_{u},d_{v}\neq1}}\sqrt{\frac{n_u+n_v-2}{n_u n_v }}\nonumber \\
&\leq& (n-5)\sqrt{\frac{n-2}{n-1}}+2\sqrt{\frac{n-3}{n-2}}+3\sqrt{\frac{1}{2}}\nonumber \\
&<& (n-5)\sqrt{\frac{n-2}{n-1}}+6\sqrt{\frac{1}{2}}. \nonumber
\end{eqnarray}
Since $G$ has maximal $ABC_{GG}$ index, it follows that $\min \{n_1,n_2\}=0$ and $\min \{m_1,m_2\}=0$. $\Box$ 

{\bf Lemma 3.5. }\  Let $G \in S_n^{r,t}$ be a connected bicyclic graph with $n \geq 6$.

(1) If $n=6$ then $  ABC_{GG}(G)\leq \sqrt{3}+\sqrt{2}+\sqrt{\frac{2}{3}}+\sqrt{\frac{4}{5}}$,
with equality holding if and only if $G= S_6^{3,3}(1,0,0)$.

(2) If $n=7$ then   $ ABC_{GG}(G)\leq\sqrt{3}+2\sqrt{2}+2\sqrt{\frac{5}{6}}$,
with equality holding if and only if $G= S_7^{3,3} (1,1,0) $.

(3) If $n=8$ then   $ ABC_{GG}(G)\leq\sqrt{2}+\sqrt{\frac{2}{3}}+\sqrt{\frac{3}{4}}+
\sqrt{\frac{4}{5}}+\sqrt{\frac{5}{12}}+3\sqrt{\frac{6}{7}}$,
with equality holding if and only if $G= S_8^{3,3}(2,1,0)$.

(4) If $n=9$ then   $ ABC_{GG}(G)\leq 2\sqrt{\frac{2}{3}}+2\sqrt{\frac{6}{15}}+
2\sqrt{\frac{4}{5}}+4\sqrt{\frac{7}{8}}$,
with  equality holding if and only if $G= S_9^{3,3} (2,2,0)$.

(5) If $10\leq n\leq15$ then   $ ABC_{GG}(G)\leq 2\left(\sqrt{\frac{2}{3}}+\sqrt{\frac{n-5}{n-4}}+
\sqrt{\frac{n-3}{3(n-4)}}\right)+(n-5)\sqrt{\frac{n-2}{n-1}}$,
with equality holding if and only if $G= S_n^{3,3}(2,2,n-9)$ .

(6) If $n \geq16$ then   $  ABC_{GG}(G)\leq 2\left(\sqrt{2}+\sqrt{\frac{n-4}{n-3}}\right)+
(n-5)\sqrt{\frac{n-2}{n-1}}$,
with equality holding if and only if $G=S_n^{3,3} (1,1,n-7)$.

{\bf Proof.}\ Let $G \in S^{3,3}_n$ be a connected bicyclic graph with maximal $ABC_{GG}$ index.
If $n=6$, then $G \in \{S_6^{3,4}(0,0,0,0,0), S_6^{3,3}(1,0,0), S_6^{3,3}(0,0,1)\}$. Since
\begin{eqnarray*}
&& ABC_{GG}(S_6^{3,3} (1,0,0))= \sqrt{3}+\sqrt{2}+\sqrt{\frac{2}{3}}+
\sqrt{\frac{4}{5}}\approx4.8572, \\
&& ABC_{GG}(S_6^{3,4}(0,0,0,0,0)) = \sqrt{3}+
2\sqrt{2}\approx4.5605, \\
&& ABC_{GG}(S_6^{3,3}(0,0,1))=2\sqrt{3}+\sqrt{\frac{4}{5}}\approx4.3585.
\end{eqnarray*}
Statement (1) follows.

If $n=7$, by  Lemma 3.1(1), 3.3, 3.4 and equation (1) we have
$$ABC_{GG}(G)\leq ABC_{GG}(S_7^{3,3}(1,1,0))=\sqrt{3}+2\sqrt{2}+2\sqrt{\frac{5}{6}}$$
with equality holding if and only if $G= S_7^{3,3}(1,1,0)$.

If $n=8$, by  Lemma 3.1(1), (2), Lemma 3.3, 3.4 and equation (1) we have
$$ABC_{GG}(G)\leq ABC_{GG}(S_8^{3,3}(2,1,0))=\sqrt{2}+\sqrt{\frac{2}{3}}+\sqrt{\frac{3}{4}}+
\sqrt{\frac{4}{5}}+\sqrt{\frac{5}{12}}+3\sqrt{\frac{6}{7}}$$
with equality holding if and only if $G= S_8^{3,3}(2,1,0)$.

If $n=9$, by Lemma 3.1(1), (2), (3), Lemma 3.3, 3.4 and equation (1) we have
$ABC_{GG}(S_9^{3,3}(3,1,0))= 2\sqrt{\frac{3}{4}}+\sqrt{\frac{3}{8}}+\sqrt{\frac{4}{5}}+
\sqrt{\frac{5}{6}}+4\sqrt{\frac{7}{8}}\approx7.8934< ABC_{GG}(S_9^{3,3} (2,2,0))=
2\sqrt{\frac{2}{3}}+2\sqrt{\frac{6}{15}}+2\sqrt{\frac{4}{5}}+4\sqrt{\frac{7}{8}}\approx8.4284$. So,
$$ABC_{GG}(G)\leq ABC_{GG}(S_9^{3,3}(2,2,0))=2\sqrt{\frac{2}{3}}+2\sqrt{\frac{6}{15}}+
2\sqrt{\frac{4}{5}}+4\sqrt{\frac{7}{8}}$$
with equality holding if and only if $G= S_9^{3,3}(2,2,0)$.

If $10\leq n\leq15$, by equation (1) and Lemma 3.1(4), Lemma 3.3, 3.4 we have
\begin{eqnarray}
ABC_{GG}(G)&\leq&ABC_{GG}(S_n^{3,3}(2,2,n-9)) \nonumber \\
&=&2\left(\sqrt{\frac{2}{3}}+\sqrt{\frac{n-5}{n-4}}+\sqrt{\frac{n-3}{3(n-4)}}\right)+
(n-5)\sqrt{\frac{n-2}{n-1}}\nonumber
\end{eqnarray}
with equality holding if and only if $G= S_n^{3,3}(2,2,n-9)$.

Finally, if $n \geq 16$, by equation (1) and Lemma 3.1(5), Lemma 3.3, 3.4 we have
$$ABC_{GG}(G)\leq ABC_{GG}(S_n^{3,3}(1,1,n-7))=2\left(\sqrt{2}+\sqrt{\frac{n-4}{n-3}}\right)+
(n-5)\sqrt{\frac{n-2}{n-1}}$$
with equality holding if and only if $G= S_n^{3,3}(1,1,n-7)$.  $\Box$

Let $C_4=v_1v_2v_3v_4v_1$ be a cycle of order 4 and $Q_4=C_4+v_1v_3$.
Denote by $B_n (n_1,n_2,n_3,n_4 )$ be the graph obtained from $Q_4$ by attaching $n_i-1 \ge 0$
isolated vertices to $v_i$ for all  $i=1,2,3,4$ with $n_1\geq n_2\geq n_3\geq n_4$. Then
 $B_{n,n-4}$ is the set of all such graphs $B_n (n_1,n_2,n_3,n_4 )$ with $n_1+n_2+n_3+n_4=5$.

{\bf Lemma 3.6. } Let $G= B_{n}(n_1,n_2,n_3,n_4)$ be a graph with $n_1, n_2, n_3, n_4 \ge 1$.
Then  $ABC_{GG}(G)< ABC_{GG}(B_{n}(n_1+1,n_2,n_3-1,n_4))$ if $n_2,n_4\ge2$ and $n_1\geq n_3 \ge 2$;
$ABC_{GG}(G)< ABC_{GG}(B_{n}(n_1,n_2+1,n_3,n_4-1))$ if $n_2 \geq n_4 \ge 2$ and $n_3=1$.

{\bf Proof}. Since all pendent edges are contribute the same to the $ABC_{GG}$ index, it
follows that
\begin{eqnarray*}
&&ABC_{GG}(B_{n}(n_1+1,n_2,n_3-1,n_4))-ABC_{GG}(G) \nonumber \\
&&=\left(\sqrt{\frac{n_1+n_2+n_4-1}{(n_1+1+n_4)n_2}}-\sqrt{\frac{n_1+n_2+n_4-2}{(n_1+n_4)n_2}}\right)\\
&&-\left(\sqrt{\frac{n_3+n_2+n_4-2}{(n_3+n_4)n_2}}-\sqrt{\frac{n_3+n_2+n_4-3}{(n_3-1+n_4)n_2}}\right)\\
&&+\left(\sqrt{\frac{n_1+n_2+n_4-1}{(n_1+n_2+1)n_4}}-\sqrt{\frac{n_1+n_2+n_4-2}{(n_1+n_2)n_4}}\right)\\
&&-\left(\sqrt{\frac{n_2+n_3+n_4-2}{(n_2+n_3)n_4}}-\sqrt{\frac{n_2+n_3+n_4-3}{(n_2+n_3-1)n_4}}\right)\\
&&+\left(\sqrt{\frac{n_1+n_3-2}{(n_1+1)(n_3-1)}}-\sqrt{\frac{n_1+n_3-2}{n_1n_3}}\right).
\end{eqnarray*}
Since $n_1\geq n_3\geq2$, it follows that
$$
\sqrt{\frac{n_1+n_3-2}{(n_1+1)(n_3-1)}}-\sqrt{\frac{n_1+n_3-2}{n_1n_3}}>0.
$$
Recalling that by Lemma 2.1 $g_{a}(x)$ is strictly increasing for $a$ if $x \ge 3$, and noticing that  $g_a(2)=0$
we deduce that
\begin{eqnarray*}
&&\left(\sqrt{\frac{n_1+n_2+n_4-1}{(n_1+1+n_4)n_2}}-\sqrt{\frac{n_1+n_2+n_4-2}{(n_1+n_4)n_2}}\right)\\
&&-\left(\sqrt{\frac{n_3+n_2+n_4-2}{(n_3+n_4)n_2}}-\sqrt{\frac{n_3+n_2+n_4-3}{(n_3-1+n_4)n_2}}\right) \\
&&=g_{n_1+n_4+1}(n_2)-g_{n_3+n_4}(n_2) \ge 0,
\end{eqnarray*}
\begin{eqnarray*}
&&\left(\sqrt{\frac{n_1+n_2+n_4-1}{(n_1+n_2+1)n_4}}-\sqrt{\frac{n_1+n_2+n_4-2}{(n_1+n_2)n_4}}\right)\\
&&-\left(\sqrt{\frac{n_2+n_3+n_4-2}{(n_2+n_3)n_4}}-\sqrt{\frac{n_2+n_3+n_4-3}{(n_2+n_3-1)n_4}}\right) \\
&&=g_{n_1+n_2+1}(n_4)-g_{n_3+n_2}(n_4) \ge 0.
\end{eqnarray*}
Therefore, the first statement is true. For the second one, we have
$$ ABC_{GG}(B_{n}(n_1,n_2+1,1,n_4-1))-ABC_{GG}(B_{n}(n_1,n_2,1,n_4)) $$
\begin{eqnarray*}
&&=\left(\sqrt{\frac{n_1+n_2+n_4-2}{(n_1+n_2+1)(n_4-1)}}-\sqrt{\frac{n_1+n_2+n_4-2}{(n_1+n_2)n_4}}\right)\\
&&-\left(\sqrt{\frac{n_1+n_2+n_4-2}{(n_1+n_4)n_2}}-\sqrt{\frac{n_1+n_2+n_4-2}{(n_1+n_4-1)(n_2+1)}}\right)\\
&&+\left(\sqrt{\frac{n_2+n_4-1}{(n_2+2)(n_4-1)}}-\sqrt{\frac{n_2+n_4-1}{(n_2+1)n_4}}\right)\\
&&+\left(\sqrt{\frac{n_2+n_4-1}{n_4(n_2+1)}}-\sqrt{\frac{n_2+n_4-1}{(n_4+1)n_2}}\right).
\end{eqnarray*}
Since $n_2 \ge n_4 \ge 2$, it follows that the first, third term are all positive, and the fourth term
is nonnegative. The second term is
also positive if $n_2 \ge n_1+n_4$, and so the second statement is true in this case. In what follows we
consider the other case when $n_2 < n_1+n_4$. It suffices to show that the sum of the first and second term is
nonnegative, or
$$
\sqrt{\frac{1}{(n_1+n_2+1)(n_4-1)}}-\sqrt{\frac{1}{(n_1+n_2)n_4}} \ge
$$
$$
\sqrt{\frac{1}{(n_1+n_4)n_2}}-\sqrt{\frac{1}{(n_1+n_4-1)(n_2+1)}}.
$$

Let $m=n_1+n_2+n_4=n-1$ and $g(x)=\frac{1}{\sqrt{x(m-x)}}-
\frac{1}{\sqrt{(x+1)(m-x-1)}}$ with $1 \le x < \frac{m}{2}$. Then $m \ge 5$,
and the above inequality becomes $g(n_4-1) \ge g(n_2)$. So, to show the
above inequality we need only show that $g(x)$ is decreasing when $1 \le x < \frac{m}{2}$.
Noticing that
$$
\frac{dg}{dx}=\frac{m-2x-2}{2\left((x+1)(m-x-1)\right)^{3/2}}-\frac{m-2x}{2\left(x(m-x)\right)^{3/2}},
$$
to show $\frac{dg}{dx}<0$ it suffices to show that
$$
(m-2x)^2(x+1)^3(m-x-1)^3>x^3(m-x)^3(m-2x-2)^2.
$$
But this is obvious since $(x+1)/(m-x)>x/(m-x-1)$ for all $x$ with $1 \le x < m/2$. And so,
the second statement is true. $\Box$ 

The above lemma tell us that in $B_{n,n-4}$, graphs with maximal $ABC_{GG}$ index are of the form
$B_n(n_1,n_2,1,1)$, where $n_1, n_2 \ge 1$. 

{\bf Lemma 3.7.} If $G=B_{n}(n_1,n_2,1,1)$, then $ABC_{GG}(G) \le 2\sqrt{\frac{n-3}{n-2}}+\sqrt{2}+
\sqrt{\frac{n-4}{n-3}}+(n-4)\sqrt{\frac{n-2}{n-1}}$, with the equality holding if and only if
$G=B_{n}(n-3,1,1,1)$. 

{\bf Proof.} Since
\begin{eqnarray*}
ABC_{GG}(G)&=&\sqrt{\frac{n-3}{(n_1+1)n_2}}+\sqrt{\frac{n_2}{n_2+1}}+
\sqrt{\frac{n_1-1}{n_1}}+\sqrt{\frac{n-3}{n-2}}\\
&&+\sqrt{\frac{1}{2}}+(n-4)\sqrt{\frac{n-2}{n-1}},
\end{eqnarray*}
\begin{eqnarray} ABC_{GG}(B_{n}(n-3,1,1,1)) &= &2\sqrt{\frac{n-3}{n-2}}+
\sqrt{\frac{n-4}{n-3}}+\sqrt{2} \nonumber \\
& &+(n-4)\sqrt{\frac{n-2}{n-1}}.\nonumber
\end{eqnarray}
we have
\begin{eqnarray}
&&ABC_{GG}(B_{n}(n-3,1,1,1))-ABC_{GG}(G)\nonumber\\
&&=\left(\sqrt{\frac{n-3}{n-2}}-\sqrt{\frac{n_2}{n_2+1}}\right)+\left(\sqrt{\frac{1}{2}}-
\sqrt{\frac{n-3}{(n_1+1)n_2}}\right)\nonumber\\
& &+\left(\sqrt{\frac{n-4}{n-3}}-\sqrt{\frac{n_1-1}{n_1}}\right).\nonumber
\end{eqnarray}
If $ n_1, n_2 \geq 2$, then $n_1, n_2 \le n-4$. And so, $ABC_{GG}(G)<ABC_{GG}(B_{n}(n-3,1,1,1))$.
This observation shows that $ABC_{GG}(G)=ABC_{GG}(B_{n}(n-3,1,1,1))$ only if $n_1=1$ or
$n_2=1$. By direct calculation the above difference one can show that $n_2=1$. And so,
the lemma follows.   $\Box$   

Let $B^{r,t}_n(s)$ be the set of all such $n$-vertex bicyclic graphs whose two cycles $C_r$ and $C_t$
have shortest common path $P_s$, and every vertex not in these two cycles is pendent,
where $s \geq 2$ and $ r \geq t \geq 3$. 

{\bf Lemma 3.8. }\ Let $G$ be a connected bicyclic graph of order $n \geq5$, whose two cycles $C_r$ and
$C_t$ have shortest common path $P_s$ with $s \ge 2$ and $r \ge t \ge 3$. Then
$ABC_{GG}(G) < ABC_{GG}(B_{n}(n-3,1,1,1))$.

{\bf Proof.}\ Assume that $G$ has maximal $ABC_{GG}$ index in all the graphs postulated in this lemma.
By the Lemma 3.2, $G \in B^{r,t}_n(s)$.

If $s \geq 3$, then $r \geq t \geq 4$, $t-s \ge 1$ and $r+t-s \ge 2 +s \ge 5$. In this case,
 $n_u, n_v \geq2$ for each $uv \in E(C_r)\cup E(C_t)$. And so,
 $$
 \sqrt{\frac{n_u+n_v-2}{n_un_v}}\leq \sqrt{\frac{1}{2}}<\sqrt{\frac{n-2}{n-1}}.
 $$
 Hence,
\begin{eqnarray}
ABC_{GG}(G)&= &\sum\limits_{\substack {uv\in E(G) \\ d_{u}=1}}\sqrt{\frac{n_u+n_v-2}{n_un_v }}+
\sum\limits_{\substack{ uv\in E(G) \\ d_{u},d_{v}\neq1}}\sqrt{\frac{n_u+n_v-2}{n_u n_v }}\nonumber \\
& \le & (n-r-t+s)\sqrt{\frac{n-2}{n-1}}+(r+t-s+1)\sqrt{\frac{1}{2}}\nonumber \\
&=&  (n-4)\sqrt{\frac{n-2}{n-1}}+ (s-r-t+4)\sqrt{\frac{n-2}{n-1}}\nonumber \\
& & +(r+t-s-4)\sqrt{\frac{1}{2}}+5\sqrt{\frac{1}{2}}\nonumber \\
&< &(n-4)\sqrt{\frac{n-2}{n-1}}+5\sqrt{\frac{1}{2}}\nonumber \\
&<&2\sqrt{\frac{n-3}{n-2}}+\sqrt{2}+\sqrt{\frac{n-4}{n-3}}+(n-4)\sqrt{\frac{n-2}{n-1}} \nonumber \\
&=& ABC_{GG}(B_{n}(n-3,1,1,1)).\nonumber 
\end{eqnarray}

If $s=2$ and $r, t \geq 4$, then the above inequalities are also true. Finally, we consider the case when
$s=2$ and $r \geq 4$, $t=3$. Assume the number of pendent vertices  adjacent to the vertex of $C_t-C_r$
is $m$. Noticing that $n \ge 5$ in this case, we have
\begin{eqnarray*}
& &ABC_{GG}(G)\leq  r \sqrt{\frac{1}{2}}+(n-r-1)\sqrt{\frac{n-2}{n-1}}+\left\{ \begin{array}{ll}
 2\sqrt{\frac{1}{2}}, & \quad \mbox{if } \ m \ge 1; \\ 2\sqrt{\frac{n-4}{n-3}}, & \quad \mbox{if}\ m=0.
\end{array}
\right. \\
& &< r\sqrt{\frac{1}{2}}+(3-r)\sqrt{\frac{n-2}{n-1}}+(n-4)\sqrt{\frac{n-2}{n-1}}+\left\{ \begin{array}{ll}
 2\sqrt{\frac{1}{2}}, & \quad \mbox{if } \ m \ge 1; \\ 2\sqrt{\frac{n-4}{n-3}}, & \quad \mbox{if}\ m=0.
\end{array}
\right. \\
& & < 3\sqrt{\frac{1}{2}}+(n-4)\sqrt{\frac{n-2}{n-1}}+\left\{ \begin{array}{ll}
 2\sqrt{\frac{1}{2}}, & \quad \mbox{if } \ m \ge 1; \\ 2\sqrt{\frac{n-4}{n-3}}, & \quad \mbox{if}\ m=0.
\end{array}
\right. \\
& & < ABC_{GG}(B_{n}(n-3,1,1,1)).
\end{eqnarray*}
The above discussion shows that $G=B_n(n_1,n_2,n_3,n_4)$. And so, the lemma follows from
Lemma 3.6 and 3.7. $\Box$

{\bf Theorem 3.9. }\ If $G$ is a connected bicyclic graph of order $n \geq 4$
then  $ABC_{GG}(G) \le 2\sqrt{\frac{n-3}{n-2}}+\sqrt{2}+\sqrt{\frac{n-4}{n-3}}+(n-4)\sqrt{\frac{n-2}{n-1}}$,
with equality holding if and only if $G=B_n(n-3,1,1,1)$.

{\bf Proof.}\ The case when $n=4$ is trivial. Suppose in what follows that $n \geq5$.
If $n=5$, by the Lemma 3.6, 3.7 and 3.8, we need only to compare $ABC_{GG}(B_5(2,1,1,1))$
and $ABC_{GG}(S_5^{3,3}(0,0,0)$. Since
$ABC_{GG}(B_5(2,1,1,1))=2\sqrt{\frac{2}{3}}+3\sqrt{\frac{1}{2}}+\sqrt{\frac{3}{4}}>
4\sqrt{\frac{2}{3}}=ABC_{GG}(S_5^{3,3}(0,0,0))$, the theorem follows in  this case. Theorem is
also true when  $n=6$ since $ABC_{GG}(B_6 (3,1,1,1))=\sqrt{3}+\sqrt{2}+\sqrt{\frac{2}{3}}+
2\sqrt{\frac{4}{5}}>\sqrt{2}+\sqrt{3}+\sqrt{\frac{2}{3}}+\sqrt{\frac{4}{5}}=
ABC_{GG}(S_6^{3,3}(1,0,0) )$.

When $n=7$, $ABC_{GG}(B_7(4,1,1,1))=2\sqrt{\frac{4}{5}}+\sqrt{2}+\sqrt{\frac{3}{4}}+
3\sqrt{\frac{5}{6}}\approx6.8077$ and $ABC_{GG}(S_7^{3,3}(1,1,0))=\sqrt{3}+2\sqrt{2}+
2\sqrt{\frac{5}{6}}\approx6.3862$; When $n=8$, $ABC_{GG}(B_8(5,1,1,1))=2\sqrt{\frac{5}{6}}+
\sqrt{2}+\sqrt{\frac{4}{5}}+
4\sqrt{\frac{6}{7}}\approx7.8377$ and $ABC_{GG}(S_8^{3,3}(2,1,0) )=\sqrt{\frac{2}{3}}+
\sqrt{\frac{3}{4}}+\sqrt{\frac{4}{5}}+\sqrt{\frac{5}{12}}+\sqrt{2}+3\sqrt{\frac{6}{7}}\approx7.4141$;
When $n=9$, $ABC_{GG}(B_9(6,1,1,1))=2\sqrt{\frac{6}{7}}+\sqrt{2}+\sqrt{\frac{5}{6}}+
5\sqrt{\frac{7}{8}}\approx8.8558$ and $ABC_{GG}(S_9^{3,3}(2,2,0) )=2\sqrt{\frac{2}{3}}+
2\sqrt{\frac{6}{15}}+2\sqrt{\frac{4}{5}}+4\sqrt{\frac{7}{8}}\approx8.4284$. And so, the theorem is also true
in all these cases.

If $10\leq n \leq15$, by Lemma 3.5, 3.6, 3.7 and 3.8, we need only to compare $ABC_{GG}(B_n(n-3,1,1,1))
=2\sqrt{\frac{n-3}{n-2}}+\sqrt{2}+\sqrt{\frac{n-4}{n-3}}+(n-4)\sqrt{\frac{n-2}{n-1}}$ and
$ABC_{GG}(S_n^{3,3}(2,2,n-9))=2 \left(\sqrt{\frac{2}{3}}+\sqrt{\frac{n-5}{n-4}}+
\sqrt{\frac{n-3}{3(n-4)}}\right)+(n-5)\sqrt{\frac{n-2}{n-1}}$. Since $n-1>n-3>n-4$,
by Lemma 2.1 we conclude that $ABC_{GG}(B_n(n-3,1,1,1))-ABC_{GG}(S_n^{3,3}(2,2,n-9))=
2\sqrt{\frac{n-3}{n-2}}-2\sqrt{\frac{2}{3}}+\sqrt{2}-2\sqrt{\frac{n-3}{3(n-4)}}+
\sqrt{\frac{n-4}{n-3}}-\sqrt{\frac{n-5}{n-4}}
+\sqrt{\frac{n-2}{n-1}}-\sqrt{\frac{n-5}{n-4}}>0$. Hence,
$$
ABC_{GG}(S_n^{3,3}(2,2,n-9))<ABC_{GG}(B_n(n-3,1,1,1)).
$$

Finally, if $n \geq16$, since $ABC_{GG}(S_n^{3,3}(1,1,n-7))=
2\left(\sqrt{2}+\sqrt{\frac{n-4}{n-3}}\right)+(n-5)\sqrt{\frac{n-2}{n-1}}$ it follows that
$ABC_{GG}(B_n(n-3,1,1,1))-ABC_{GG}(S_n^{3,3}(1,1,n-7))=2\sqrt{\frac{n-3}{n-2}}-
\sqrt{2}+\sqrt{\frac{n-2}{n-1}}-\sqrt{\frac{n-4}{n-3}}>0$. Hence,
$$
ABC_{GG}(S_n^{3,3}(1,1,n-7))<ABC_{GG}(B_n(n-3,1,1,1)).
$$
By Lemma 3.5, 3.6, 3.7 and 3.8, the theorem follows.  $\Box$

\end{document}